\documentclass [12pt] {article}

\begin{document} 
\newtheorem{example}{Example}[section] 
\newtheorem{theorem}{Theorem}[section]
\newtheorem{lemma}{Lemma}[section]
\newtheorem{corollary}{Corollary}[section]
\newtheorem{proposition}{Proposition}[section]
\newtheorem{remark}{Remark}[section]
\setlength{\textwidth}{13cm}
\def \halmos{\hfill\mbox{ qed}\\}
\newcommand{\eqnsection}
{\renewcommand{\theequation}{\thesection.\arabic{equation}}
\makeatletter \csname  @addtoreset\endcsname{equation}{section}
\makeatother}

%
%\newtheorem{lemma}[theorem]{Lemma}
%\newtheorem{remark}[theorem]{Remark}
% \newtheorem{example}[theorem]{Example}
% \newtheorem{corollary}[theorem]{Corollary}
%\newtheorem{conjecture}[theorem]{Conjecture}
%\newcommand{\eqnsection}
%{\renewcommand{\theequation}{\thesection.\arabic{equation}}
%\makeatletter \csname  @addtoreset\endcsname{equation}{section}
%\makeatother}

\def \njc{{\bf !!! note Jay's change  !!!}}
\def  \enc{{\bf end of current changes }}
\def \nnc{{\bf !!! note NEW	 change  !!!}}
\def\square{{\vcenter{\vbox{\hrule height.3pt
                    \hbox{\vrule width.3pt height5pt \kern5pt
                       \vrule width.3pt}
                    \hrule height.3pt}}}}
\def \grad{\bigtriangledown}
\def \nc{{\bf !!! note change !!! }}
        \def \Proof{\noindent{\bf Proof}$\quad$}
% Special Commands
\newcommand{\re}[1]{(\ref{#1})}
\def \ov{\overline}
\def \un{\underline}
\def \be{\begin{equation}}
\def \ee{\end{equation}}
\def \bt{\begin{theorem}}
\def \et{\end{theorem}}
\def \bc{\begin{corollary}}
\def \ec{\end{corollary}}
\def \br{\begin{remark} }
\def \er{ \end{remark}}
\def \bl{\begin{lemma}}
\def \el{\end{lemma}}
\def \bex{\begin{example}}
\def \eex{\end{example}}
\def \bea{\begin{eqnarray}}
\def \eea{\end{eqnarray}}
\def \bas{\begin{eqnarray*}}
\def \eas{\end{eqnarray*}}
% Greek Letters
\def \al{\alpha}
\def \bb{\beta}
\def \ga{\gamma}
\def \Ga{\Gamma}
\def \de{\delta}
\def \De{\Delta}
\def \ep{\epsilon}
\def \vep{\varepsilon}
\def \la{\lambda}
\def \La{\Lambda}
\def \ka{\kappa}
\def \om{\omega}
\def \Om{\Omega}
\def \va{\varrho}
\def \ffi{\Phi}
\def \vf{\varphi}
\def \si{\sigma}
\def \Si{\Sigma}
\def \vsi{\varsigma}
\def \th{\theta}
\def \Th{\Theta}
\def \ups{\Upsilon}
\def \ze{\zeta}
\def \tr{\nabla}
% Math Symbols
\def \ff{\infty}
\def \wh{\widehat}
\def \wt{\widetilde}
\def \dar{\downarrow}
\def \rar{\rightarrow}
\def \uar{\uparrow}
\def \sbs{\subseteq}
\def \mpt{\mapsto}
\def \R{{\bf R}}
\def \G{{\bf G}}
\def \H{{\bf H}}
\def \Z{{\bf Z}}
\def \S{{\bf
S}}
\def \sfB{{\sf B}}
\def \sfS{{\sf S}}
\def \T{{\bf T}}
\def
\C{{\bf C}}
\def \AA{{\mathcal A}}
\def \BB{{\mathcal B}}
\def
\CC{{\mathcal C}}
\def \DD{{\mathcal D}}
\def \EE{{\mathcal E}}
\def
\FF{{\mathcal F}}
\def \GG{{\mathcal G}}
\def \HH{{\mathcal H}}
\def
\II{{\mathcal I}}
\def \JJ{{\mathcal J}}
\def \KK{{\mathcal K}}
\def
\LL{{\mathcal L}}
\def \MM{{\mathcal M}}
\def \NN{{\mathcal N}}
\def
\OO{{\mathcal O}}
\def \PP{{\mathcal P}}
\def \QQ{{\mathcal Q}}
\def
\RR{{\mathcal R}}
\def \SS{{\mathcal S}}
\def \TT{{\mathcal T}}
\def
\UU{{\mathcal U}}
\def \VV{{\mathcal V}}
\def \ZZ{{\mathcal Z}}
\def \Pxh{P^{x/h}}
\def \Exh{E^{x/h}}
\def \Px{P^{x}}
\def \Ex{E^{x}}
\def
\Prh{P^{\rho/h}}
\def \Erh{E^{\rho/h}}
\def \p{p_{t}(x,y)}
\def
\({\left(}
\def \){\right)}
\def \lk{\left[}
\def \rk{\right]}
\def
\lc{\left\{}
\def \rc{\right\}}
\def \bsq{\ $\Box$}
\def
\nn{\nonumber}
\def \Bo{\bigotimes}
\def \bo{\times}
\def
\ot{\times}
\def
\bs{\begin{slide} }
\def \es{\end{slide} }
\def \bpr{\begin{proof}
}
\def \epr{\end{proof} }
\def \cd{\,\cdot\,}
\def
\st{\stackrel{def}{=}}
\def \as{almost surely }
\def \fix {{\bf
!!!!! }}
\def \stl{\stackrel{law}{=}}
\def \std{\stackrel{dist}{=}}
\def \stdto{\stackrel{dist}{\longrightarrow}}
\def  \enc{{ \bf end of current changes}}
%%%%%%%%%%%%%%%%%
\def\square{{\vcenter{\vbox{\hrule height.3pt
                    \hbox{\vrule width.3pt height5pt \kern5pt
                       \vrule width.3pt}
                    \hrule height.3pt}}}}
\def\qed{{\hfill $\square$ }}

           \def \tb{|\!|\!|}

        \eqnsection
\bibliographystyle{amsplain}

\title{  An almost sure limit theorem for Wick powers of Gaussian differences quotients}  
        \author{ Michael B. Marcus\,\, Jay Rosen \thanks{Research of 
both authors
supported by  grants from the National Science Foundation and PSCCUNY.}}

%\date{ }

\maketitle
\eqnsection

\bibliographystyle{amsplain}

          \def \tb{|\!|\!|}

\begin{abstract}   
LetÊ $G=\{G(x),x\in R_+ \}$,Ê Ê $G(0)=0$, be a mean zeroÊ GaussianÊ process with $E(G(x)-G(y))^2=\si ^2( x-y )  $. Let $
 \rho (x)=     \frac12{d^{2}\over dx^2}\si^2(x)$, $x\ne 0 $. When $\rho^{k}$ is integrable at zero and satisfies some additional regularity conditions,
 \[ \lim_{h\downarrow 0}Ê  \int :\(\frac{G(x+h)-G(x)}{h}
\)^{k }:g(x)\,dx
Ê =\,:(G') ^{ k}:( g)\hspace{.3 in}a.s.  
\]
 for all $g\in \BB_{0}(R^{+})$, the set of bounded Lebesgue measurable functions on $R_+$ with compact support. Here $G'$  is a generalized derivative of $G$ and $:(\cd)^{k}:$ is the  $k$--th order Wick power.
\end{abstract}

\section{Introduction} LetÊ $G=\{G(x),x\in R_+ \}$,Ê Ê $G(0)=0$, be a mean zeroÊ GaussianÊ process
with stationary increments, and set 
\be 
Ê Ê ÊÊ E(G(x)-G(y))^2=\si ^2( x-y ) \label{m1.2}=\si
^2(|x-y|). \label{1.1}
\ee 
 (The function $\si^2$ is referred to as the
increment's variance of $G$.)Ê Ê We assume that
 \be
\mbox{$\si^2(h)$ is a convex function that is regularly varying at zero;}\label{rv}
\ee
\be\lim_{h\to 0}{h^2\over \si^2(h)}= 0\qquad \mbox{and}\qquad \lim_{h\to
0}{
\si^2(h)\over h}= 0;\label{mm5.17o}
\ee
\be
{ \si^2(s+h)+ \si^2(s-h)-2\si^2(s)\over h^2}\le C{\si^2(s)\over
s^2}\quad\mbox{for } h\le { s\over 8 };\label{mm5.9}
\eeÊ
 \be
\mbox{$\si^2(s)$ has a second derivative for each $s\ne	0$.}\label{rv2}
\ee
Note that by  (\ref{rv})  
\be \rho (s):=\frac12{d^{2}\over ds^2}\si^2(s)\geq 0.\label{rho}
\eeÊ

  \medskip	
It follows from the second condition in (\ref{mm5.17o}) that $G$ has a continuous version; (see \cite[Lemma 6.4.6]{book}). We work with this version. 
 However, when  $\lim_{x\to 0}\rho(x)=\ff$, $G $ is not   differentiable; it is not even mean square differentiable.  It is a natural question to ask whether the weak limit   
\begin{equation}
\lim_{h\rar 0}\int \({G(x+h) -G(x)\over h}\)g(x)\,dx\label{wl.1}
\end{equation}
exists in some sense. 
Here  $g\in \BB_{0}(R_{+})$, the set of bounded Lebesgue measurable functions on $R_+$ with compact support. 

We    show in \cite[Theorem 2.1]{ae} that when  $G$ satisfies the second condition in (\ref{mm5.17o})  there exists a mean zero Gaussian field $\{G'(g), g\in\BB_{0}(R_{+})\}$    with covariance   
\begin{equation}
 E\(G'(g)G'(\wt g)\)=\int\int\rho (t-s)\,g( s)\,\wt g( t)   \,ds\,dt\label{new.5}
\end{equation}
such that
\begin{equation}
\lim_{h\rar 0}\int \({G( x+h)- G( x) \over h}\)\,g( x)\,dx=   G'(g)\hspace{.15 in}\mbox{ in }
L^{2}.\label{si.79}
\end{equation}
 Because of this we  think of $G'$ as a generalized derivative of $G$.

More generally, one may  consider  
\begin{equation}
\lim_{h\rar 0}\int \({G(x+h) -G(x)\over h}\)^{k}g(x)\,dx\label{wl.2}
\end{equation}
for  any integer $k\geq 1$. 
However,   when $k$ is even, the expectation of the square of the integral in (\ref{wl.2}) contains terms in $\si^{2}(h)/h^{2}$ which goes to infinity as $h$ goes to zero by   (\ref{mm5.17o}).     To obtain a finite limit in (\ref{wl.2})     we replace $  \displaystyle  \({G(x+h) -G(x)\over h}\)^{k}  $ by a $k$--th order polynomial  \be
 \sum_{j=0}^{k}a_{j}(h)\({G(x+h) -G(x)\over h}\)^{j}
\ee where,  $a_{j}(h)$ is a non--random function of $h$,  which,   necessarily,  has the property that, at least for some $0\le j\le k$,   
$\lim_{h\rar 0}|a_{j}(h)|= \ff$. We call this process renormalization. The renormalization we use is known as the  $k$--th Wick power.

The $k$--th Wick power of a mean zero Gaussian random variable $X$ is
\be
:X^k:\,=\sum_{ j=0}^{[ k/2]}( -1)^{ j}{Ê k\choose 2 j}E( X^{ 2 j})\,\,X^{ÊÊ k-2j
}.\label{defin}
\ee
When $X=N(0,1)$,Ê  
\be 
 :X^{ k}:\,=\sqrt{k!\,}H_{k}( X),\label{1.14}
 \ee
  where $H_{k}$ is the $k$--th Hermite polynomial.    One advantage of  Wick powers over  Hermite polynomials is that they are homogeneous, i.e., for $a\in R^{1}$,
\begin{equation}
 :(aX)^k:=a^{k}:X^k:.  \label{1.15}
   \end{equation} Therefore,Ê when
 $X$ has variance  $\si^2_X$,  
\be
:X^k:\,=\sqrt{k!\,}\si_X^kH_{k}\( {X\over \si _X}\).\label{0.9}
\ee

  When $\rho^{k}$ is locally integrable and bounded away from the origin we   construct   a $k$--th order Wick power Gaussian chaos  from the mean zero Gaussian field   $G' =\{G'(f ),f\in\BB_{0}(R_{+})\}$ in the following way: 
  For each $\de\in (0,\de_{0}]$, for some $\de_{0}>0$, let $f_{\de}(s)$ \label{page3} be a continuous positive symmetric function on
$(s,\de)\in R_+ \times (0,1]$, with support in the ball of radius $\de$, with
$\int f_{\de}(y)\,dy=1$. That is,
$f_{\de}$ is a   continuous  approximate identity.       In \cite[( 3.25) and (3.26)]{ae} we show that  for the Gaussian processes 
$G$ considered here,    for all $g\in\BB_{0}(R_{+})$,
\be 
:(G')^k:(g):=\lim_{\de\rar 0}\int
:(G' ( f_{x,\de}))^k: g(x)\,dx\quad \mbox{  in $L^2$  }\label{wwc}
\ee
and
\be
E\( :(G')^k:(g) \)^2=k!\int\!\!\int \rho^{k}(x-y) g(x)g(y)\,dx\,dy\label{3.7b}. 
\ee  
In \cite[Theorem 3.1]{ae} we show that
\begin{equation}
  \lim_{h\to 0}  \int
:\({G( x+h)- G( x) \over h}\)^{ k}:\,g( x)\,dx= \,\, :(G')^{ k} :( g) \qquad \label{nor4.6}
\mbox{in
$  L^2$}.
\ee

(For example, if $\si^{2}(h)=h^{r}$, in order for $\rho^{k}$ to be  locally integrable and to have $\lim_{x\to 0}\rho(x)=\ff$ as required by the first condition in (\ref{mm5.17o}),  it is necessary that $  \displaystyle  \frac{2k-1}{k}<r<2$.)

\medskip	In this paper we obtain the rather remarkable result that, under some additional mild regularity conditions on $\rho$, the limit in (\ref{nor4.6}) is almost sure.
  
   \bt \label{theo-9.2} LetÊ
$G=\{G(x),x\in R_+\}$, $G(0)=0$  be a mean zero Gaussian processÊ with
stationary incrementsÊ satisfying   (\ref{1.1})--(\ref{rho}). Fix an integer $k\ge 1$ and   assume that Ê  
there exists a
$ \,0<\de<1/2 $ and an
$M>0$ such that
\begin{equation}
\rho( x)\leq {C_{M} \over |x|^{(1 -Ê \de)/k }}:=C_M\,\varphi(|x|),\quad\hspace{.2 in} 0<|x|\leq
M
\label{nor14.5a}
\end{equation}
and 
\begin{equation}
 |\rho( x+h) -\rho( x)|\leq C_{M}{|h|Ê \over |x| }\,\rho( x),
\hspace{.2 in}4|h|\leq |x|\leq M. \label{nor14.5}
\end{equation}
 Then Êfor all $g\in \BB_{0}(R_{+})$,  
\begin{equation} \lim_{h\downarrow 0}Ê  \int :\(\frac{G(x+h)-G(x)}{h}
\)^{k }:g(x)\,dx
Ê =\,:(G') ^{ k}:( g)\hspace{.3 in}a.s. \label{nor14.6}
\end{equation}
\et
      
    (Note that by using Wick powers it is clear that (\ref{nor14.6}) deals with a generalized derivative of $G$. This point would be obscured if we expressed (\ref{nor14.6}) in terms of Hermite polynomials.)

  \medskip	 For a fixed $g\in \BB_{0}(R_{+})$ both the left-hand side and right-hand side of     (\ref{nor14.6})
    are $k$--th order Gaussian chaoses.  Let   $\{:X_{h}^{k}:(g), h\in (0,1]\}$, denote the left-hand side of     (\ref{nor14.6}) and $:X_{0}^{k}:(g)$ denote the right-hand side of     (\ref{nor14.6}). Theorem \ref{theo-9.2} is the statement that for all $g\in\BB_{0}(R_{+})$, the   $k$-th order Gaussian chaos process
\be
{\cal{X}}: =\{ X_{h} , h\in [0,1]\}:=\{: X_{h}^{ k}: ( g), h\in [0,1]\},\label{chi}
\ee
 has  a continuous version. Of course there is no problem in choosing the version. The process on the left in (\ref{nor14.6}) is continuous in $h\in (0, 1]$. We can take $:(G') ^{ k}:( g)$ to be its limit on the set of probability one for which the limit exists, and to be zero otherwise.

\medskip	  Theorem \ref{theo-9.2} is proved in Section \ref{sec6} using a majorizing measure result for the continuity of Gaussian chaoses.  Technically, this is an interesting application of this theory, because  the proof consists of  obtaining  continuity
at a single point. To prove (\ref{nor14.6}) we need a
majorizing measure condition for exponential Orlicz spaces based on the function $\exp x^{q}-1$ for $q\le 1$.  Whereas it is   known that such results exist we could not find a
reference, so we provide proofs in Section \ref{orlicz}.Ê We prove Theorem
\ref{theo-9.2} in Section
\ref{sec6}
 
 \section{$L^{2}$ results}
 
 We list here some $L^{2}$ estimates   we need in this paper that are obtained in \cite{ae}. To better motivate these results we state the main result in \cite{ae} and explain how it led to our consideration of
 Theorem \ref{theo-9.2}   in this paper. 
 
 \bt[Theorem 1.1, \cite{ae}]  \label{theo-2.1} Let  $f$ be a    function with 
$Ef^{2}(\eta)<\ff$, where $\eta=N(0,1)$. Then under the hypotheses of Theorem \ref{theo-9.2}
\bea
&& \int_a^bf\(\frac{G(x+h)-G(x)}{\si (h)}
\)\,dx \label{abst}
 \\
&&\qquad  = \sum_{j=0}^{k} (h/\si(h))^{j}\,\,{E(H_{j}(\eta)
f(\eta))\over\sqrt {j!}}\,\,:(G')^{j}:(I_{[a,b]})\, +o\({h\over\si (h)}\)^{k}\nn
\eea
in  $L^2$. Here $H_j$ is the $j$-th Hermite
polynomial    and 
  $:(G')^{j }:(I_{[a,b]})$ is a $j $-th order Wick power Gaussian
chaos as described in (\ref{wwc}).
\et  
 
We wondered whether (\ref{abst}) could be almost sure. In Theorem \ref{theo-9.2} we show that when $f(\cd)=H_{k}(\cd) $ it is. Note that in this case the right-hand side of (\ref{abst}) is  
\begin{equation}
    {(h/\si(h))^{k}\over \sqrt {k!}}\,\,:(G')^{k}:(I_{[a,b]})\, +o\({h\over\si (h)}\)^{k}.
   \end{equation}
   and  by (\ref{1.14}) and (\ref{1.15}) the left-hand  side of (\ref{abst}) is  
   \begin{equation}
   Ê     {(h/\si(h))^{k}\over \sqrt {k!}} \int_{a}^{b} :\(\frac{G(x+h)-G(x)}{h}
\)^{k }: \,dx
   \end{equation}
Thus Theorem \ref{theo-2.1} gives the limit in $L^{2}$ of  (\ref{nor14.6}) when $g=I_{[a,b]}$.  
 
 \medskip	The next lemma which is  part of  \cite[Lemma  4.2]{ae} provides part of the $L^{2}$ metric estimates that are needed in proof of continuity of $\cal{X} $.
 
  \bl \label{calc}Let  
$G=\{G(x),x\in R_+\}$, $G(0)=0$, be a mean zero Gaussian process  with
stationary increments  and set $\si ^2(|x-y|)=  E(G(x)-G(y))^2$.  Set $\rho (s)=
\frac12{d^{2}\over ds^2}\si^2(s) $.  Fix  an integer $j_0\ge 1$ and assume that  
there exists a  
$ \,0<\de<1$ and an
$M>0$ such that  (\ref{nor14.5}) holds, and (\ref{nor14.5a}) holds with $k$ replaced by $j_{0}$.
Then  for $1\le j\le j_0$   and any   $g\in\BB_{0}(R_{+})$,
\be
\|:X_{h}^{ j}:( g)-:X_{0}^{ j}:( g)\|_{2}\le C(|h|\varphi^j(h))^{1/2}.\label{.03}
\ee
\el

\section {Continuity conditions for stochastic processes in 
exponential Orlicz spaces}
 \label{orlicz}

Let $\|\cd\|_{\psi_q}$ÊÊ denote the norm in the
Orlicz space $L^{\psi_q}(dP)$, where
\be
\label{psi}
\psi_q(x)= \left\{ \begin{array}{l@{\quad\quad}r}
\exp(x^q) -1Ê & 1\le q<\ff \\
\exp\exp(x) -eÊ & q=\ff. \end{array} \right.\quad x\in R^+
\ee
    For
$0<q<1$, we define
\be
\label{psi-1}
\psi_q(x)= \left\{ \begin{array}{l@{\quad\quad}r}
K_q\,xÊ &0\leÊ x<\( \frac{1}{q}\)^{1/q}
\\
\exp(x^q) -1Ê & x\ge\( \frac{1}{q}\)^{1/q}Ê \end{array} \right.
\ee
where
\be K_q={\exp( x_0 ^q) -1 \over x_0}\qquad \mbox{and
$\qquad x_0:=x_0(q)=(1/q)^{1/q}$},
\ee
so that $\psi_q(x)$ is continuous.

\bl
For
$0< q< \ff$, $\psi_q(x)$ isÊÊ convex and increasing
andÊ there exists a constant
$C_q<\ff$,ÊÊ for which
\be
\psi_q(x)\le C_q\( \exp(x^q) -1\)\label{psiq}
\ee
 and
\be
\exp(x^q)\le C_q\(\psi_q(x)Ê +1\).\label{psiqj}
\ee
 In addition $C_q=1$ for $1\le q<\ff$.
\el

\ProofÊ This is trivial when $1\le q< \ff$. We consider the other cases. To
show $\psi_q(x)$ isÊÊ convex we show that its derivative is increasing. It is
easy to check that the derivative of $\psi_q(x)$ from the leftÊ at $x_0$ is
less than the derivative from the right at $x_0$. It is also easy to check that
the second derivative
 of $\psi_q(x)$ is
positive for $x\ge ((1-q)/q)^{1/q}$. Therefore, the derivative of $\psi_q(x)$
is increasing on $[Ê x_0,\ff)$, so $\psi_q(x)$ isÊÊ convex.

Since $\exp(x^q) -1\geq x^q$ for all $x\geq 0$, we see that
(\ref{psiq}) holds for $0\leq x\leq 1$ with $C_q=K_q$. Similarly, choosing $m$
so that $mq\geq 1$, $\exp(x^q) -1\geq x^{mq}/m!$, so that
(\ref{psiq}) holds for $1\leq x\leq x_{0}$ with $C_q=m!K_q$. It is then clear thatÊ
 (\ref{psiq})Ê with $C_q=\max (m!K_q,1)$ holds for all $x$ .Ê By
further increasing $C_q$ÊÊ it is easy to see that (\ref{psiqj}) also
holds.\qed

\medskipÊ We note the following obvious relationships:
\bl\label{inverses}
\be
\label{2mm.psi}\psi^{-1}_q(x)=\left\{
 \begin{array}{l@{\quad\qquad}l}
 \log\log(e+x) & q=\ffÊ Ê \\
\(\log(1+x)\) ^{1/q}&Ê 1\le q<\ff,
\end{array} \right.
\ee
andÊ for $0<q<1$
\be
\label{3mm.psi}\psi^{-1}_q(x)=\left\{
 \begin{array}{l@{\quad\qquad}l}
x/K_q &0\le x\le K_qx_0 \\
\(\log(1+x)\) ^{1/q}&Ê K_qx_0<x\le\ff.
\end{array} \right.
\ee
 \el

\medskip
For each $0<q\le\ff$ letÊ $L^{\psi_q}(\Om,P)$ denote the set ofÊ random
variables
$\xi:\Om\to C$ such that $E\psi_q\( |\xi|/c \)<\ff$ for some $c>0$.
$L^{\psi_q}(\Om,P)$ is a Banach space with norm given by
\be
\|\xi\|_{\psi_q}=\inf\left\{c>0:E\psi_q\( |\xi|/c\)\le1\right\}.
\ee

Let
$(T,d)$ be a pseudometric space. We use $B_{d}(t,u)$, or simply
$B(t,u)$, to denote a closed ball of radius $u$ in $(T,d)$.

\bt\label{maj}ÊÊ Let $X= \{X(t): t\in T\}$Ê be a measurable separable
stochastic process on aÊ separable metric or pseudometric space
$(T,d)$ with finite diameter $D$. Suppose thatÊ $X(t)
\in L^{\psi_q}(\Om,P)$ and $\|X(t)-X(s)\|_{\psi_q} \le d(t,s)$ for all $s,t \in
T$.ÊÊ Let $0<q<\ff$ and suppose also that there exists
a probability measure $\mu$Ê Ê on $T$ such that\be
\sup_{t\in T}\int_0^D\(\log\frac1{\mu(B(t,u))}\)^{1/q}\,du<\ff.
\label{3.6.88}
\ee
Then there exists a version $X'=\{X'(t),t\in T\}$ of $X$ such that
\be
E\sup_{t\in T}X'(t)\leÊ C\sup_{t\in
T}\int_0^D\(\log\frac1{\mu(B(t,u))}\)^{1/q}\,du
\label{3.6.88b}
 \ee
 for some $C<\ff$.
Furthermore, if
\be
\lim_{\ep\to 0}\sup_{t\in
T}\int_0^\ep\(\log\frac1{\mu(B(t,u))}\)^{1/q}\,du=0 \label{3.6.88a},
\ee then $X'$ is uniformly continuous on $T$ almost surely and
there exists aÊ positive
random variable
$Z \in L^{\psi_q}(\Om,P)$ such that  
\be
Ê ÊÊ \sup_{\stackrel{s,t\in T}{d(s,t)\le \de}} |X' ( s,\om)-X'( t,\om)|
Ê \le Z(\om) \,\sup_{s\in T} \int_0^{\de}
\(\log\frac{1}{\mu(B(s,u))}\)^{1/q} \,du\label{3.6.90h}.
\ee
almost surely.
When $q=\ff$ the results continue to hold when the above integrands are
replaced by
\be
\log^+\log\(Ê \frac1{\mu(B(t,u))}\).
\ee
\et

(The statementÊ $Z \in L^{\psi_q}(\Om,P)$ means thatÊ $\|Z\|_{\psi_q} \le
K_q$, a constant depending only on $q$.)

 \medskip We get the following useful corollary of Theorem \ref{maj}
\bc Under the hypotheses of Theorem \ref{maj} there exists a constant $C_q$
for which
\begin{equation}
\| \sup_{\stackrel{s,t\in T}{d(s,t)\le \de}} |X' ( s )-X'( t )|\,\,Ê \|_{\psi_q}\le
C_q \sup_{s\in T} \int_0^{\de}
\(\log\frac{1}{\mu(B(s,u))}\)^{1/q} \,du\label{later2}
\end{equation}
and. for any $t_0\in T$,
\begin{equation}
\| \sup_{s\in T} |X' ( s ) |\,\,Ê \|_{\psi_q}\le \|Ê X' ( t_0 ) Ê
\|_{\psi_q}+ C_q \sup_{s\in T} \int_0^{D}
\(\log\frac{1}{\mu(B(s,u))}\)^{1/q} \,du\label{later1}.
\end{equation}
\ec

\Proof The statement in (\ref{later2}) follows immediately from
(\ref{3.6.90h}). The statement in (\ref{later1}) follows from
(\ref{3.6.90h}) by writing
\bea
 \sup_{s\in T} |X' ( s ) |&\le& \sup_{s\in T} |X' ( s )-X'(t_0) |+|X'(t_0) |\\
&\le& \sup_{s,t\in T} |X' ( s )-X'(t) |+|X'(t_0) |,\nn
\eea
and using the triangle inequality with respect to $\| \cd
Ê \|_{\psi_q}$.\qed

\medskip The hypotheses of Theorem \ref{maj} are satisfied by Gaussian
processes when
$q=2$.Ê In this case it contains ideas which originated in an important early
paper byÊ Garcia,Ê Rodemich andÊÊ Rumsey Jr., \cite{GRR} andÊÊ were
developed further byÊÊ Preston, \cite{Preston1, Preston2} andÊ Fernique,
\cite{Fernique1}.  The fact that it can be extended to processes in
exponential Ê Orlicz spaces for $1\le q\le\ff$ is, no doubt, understood by many
researchers in the field of probability on Banach spaces. For lack of a
suitable reference  a proof was given in 
\cite{Marcus-Rosinski}.   

In this paper we need an extension to  
$\psi_q(x)$ for $0<q\le\ff$.  Here too we're sure many
researchers are aware that this can be done, but, once again, we have
no reference. When
$0<q<1$,  
$\exp(x^q)-1$ is not convex, so a bit more care is necessary. The key point is
the following lemma:

\bl
\label{a2} For $0< q\le \ff$,ÊÊ let $X= \{X(t): t\in T\}$ \ be a measurable
separable stochastic process on a precompact metric space $(T,d)$ such that
$\|X(t)\|_{\psi_q}\le 1$ for all $t \in T$.Ê Then there exists a random
variable
$Z$ with $\|Z\|_{\psi_q} \le C'_q$, such that for every
probability measure $m$ on $T$ andÊÊ function
$h: T\mapsto \R_+$ with $\int_{T} h(v)\, m(dv)Ê <\ff$
\bea
\label{hX}
&&\int_{T} |X(t)|h(t) \, m(dt) \\
&& \qquad \le Z
\int_{T} h(t) \Phi_q \(\(\int_{T} h(v)\, m(dv)\)^{-1} h(t)\)
\, m(dt)\nn
\eea
where
\be
\label{4mm.psi}\Phi_q=\left\{
 \begin{array}{l@{\quad\qquad}l}
 \log\log(e+x) & q=\ffÊ Ê \\
\(\log(1+x)\) ^{1/q}&Ê 1\le q<\ffÊ \\
\(2\log(1+(x/G_q))\) ^{1/q}&Ê 0\le q<1,
\end{array} \right.
\ee
where $G_q>0$.
\el

\Proof
Define
\be
\wt Z(\om)= \inf\{\al>0: \int_{T} \psi_q(\al^{-1}|X(t)|) \, m(dt) \le
1\}.\label{bsf}
\ee
We first show that
\be
\label{bZ}
\|\wt Z\|_{\psi_{q}} \leÊ C'_{q}< \ffÊ \qquadÊ 0< q \le \ff.
\ee
Let $0< q<\ff$; then for $u\ge 1$,
\begin{equation}
P\( \wt Z > u \) \leÊ P\( \int_{T} \psi_q(u^{-1}|X(t)|) \, m(dt) > 1\)  
\end{equation}
and by (\ref{psiq})
\bas
&&
P\( \int_{T} \psi_q(u^{-1}|X(t)|) \, m(dt) > 1\)Ê \\ &&\qquad \le P\(C_q \int_{T}
\exp(u^{-q}|X(t)|^q) \, m(dt) > 1+C_q\) \\
 &&\qquad =Ê P\( \(C_q\int_{T} \exp(u^{-q}|X(t)|^q) \,
m(dt)\)^{u^q} > ( 1+C_q)^{u^q}\) \\ &&\qquad \leÊ P\( C_q^{u^q}\int_{T}
\exp(|X(t)|^q) \, m(dt)Ê >Ê (1+C_q)^{u^q}\) \\Ê &&\qquad \le
\({1+C_q\overÊÊ C_q}\)^{-u^q} E\int_{T}Ê \exp(|X(t)|^q) \, m(dt)ÊÊ \\ÊÊ &&\qquad \le
\({1+C_q\overÊÊ C_q}\)^{-u^q}C_q\(E\int_{T} \psi_q(|X(t)|)Ê \,
m(dt) +1\)
\\ &&\qquad \le 2C_q
\({1+C_q\overÊ C_q}\)^{-u^q} .
\eas
Ê The fourth line follows from Jensen's inequality, the sixthÊ from (\ref{psiqj}),
and the last because
$\|X(t)\|_{\psi_q}\le 1$. Thus we get \re{bZ} when $0\le q<\ff$.

Now let $q=\ff$. Note that for each $u\ge 1$, the function $\phi_{u}(x)=\exp((\log x)^{u})$ is
 convex forÊ $x\ge e$. Using Jensen's inequality again we get that for $u\ge 1$
\bea
P\{ \wt Z > u \} &\le & P\{ \int_{T} \exp(\exp(u^{-1}|X(t)|)) \, m(dt) > e+1\} \\
&= &Ê P\{ \phi_{u}\Big(\int_{T} \exp(\exp(u^{-1}|X(t)|)) \, m(dt)\Big) >
 \phi_{u}(e+1)\} \nn\\
&\le &Ê P\{ \int_{T} \exp(\exp(|X(t)|)) \, m(dt) > \exp(\exp(cu))\} \nn\\
&\le & \exp(-\exp(cu)) \int_{T} E \exp(\exp(|X(t)|)) \, m(dt)ÊÊ \nn\\
&\le & (1+e)\exp(-\exp(cu))\nn
\eea
where $c= \log(\log(e+1))>0$. Thus we get \re{bZ} when $q=\ff$.

We now prove \re{hX}. For $1\le q\le\ff$ we have
\be
\label{Yo}
xy \le \psi_q(x) + y \psi_q^{-1}(y) \quad x, y \ge 0.
\ee
To obtain (\ref{Yo}) we first note that $\psi'_q(x)\geÊ \psi_q(x)$.  To see this
set
$h(x):=\psi'_q(x)-\psi_q(x)$. We get the desired inequality because  
$h(0)$=0 and
$h'(x)>0$ for all $x\in R_+$. To prove this last point it suffices to show that 
\be
g(x):=\frac{q-1}{x}+qx^{q-1}-1\ge 0 .\label{92.3}
\ee
To verify  (\ref{92.3}) note that the minimum of $g(x)$ takes place at
$x_1=(1 /q)^{1/q}$ and $g(x_1)>0$.

The inequality in  (\ref{Yo})ÊÊ follows from Young's inequality since
$ \psi_q(x)$ is convex and $\psi'_q(x)\geÊ \psi_q(x)$. (Recall that the final
term in (\ref{Yo}) can be taken to be
$\int_0^y (\psi'_q)^{-1}(y)\,dy$. Since $\psi'_q(x)\geÊ \psi_q(x)$,
$(\psi'_q)^{-1}(y)\leÊ (\psi_q)^{-1}(y)$, and since $(\psi_q)^{-1}(y)$ is
increasing we get (\ref{Yo}).)

When $0<q<1$ it follows from Lemma \ref{lem-5.4},   which is given at the end of this section, that
\be
\label{Yoq2}
xy \le \psi_q(x) +Ê y\(2\log (1+y/G_q)\)^{1/q}ÊÊ \quad x, y \ge 0,
\ee
for some constantÊ $G_q>0$. Therefore it follows from (\ref{Yo}),
(\ref{Yoq2}) and Lemma \ref{inverses} that
\be
\label{Yoq2w}
xy \le \psi_q(x) +Ê y\Phi_q(y)ÊÊ \quad x, y \ge 0.
\ee

Let $h: T \mapsto R_+$ be as in
the lemma. Putting $x=\wt Z^{-1}|X(t)|$ and $y=(\int_{T} h(v)\,
m(dv))^{-1}h(t)$ inÊÊ \re{Yoq2w} we get
\bea
|X(t)|h(t) &\le& \wt Z\int_{T} h(v)\, m(dv) \psi_q(\wt Z^{-1}|X(t)|) \\
&& \qquad +\wt Z h(t)\Phi_q((\int_{T} h(v)\, m(dv))^{-1} h(t)).\nn
\eea
Integration with respect to $m$, and using the definition (\ref{bsf}), gives
\bea
&&\int_{T} |X(t)|h(t) \, m(dt)\label{eee} \\
&&\qquad\leÊ \wt Z\int_{T} h(t)\, m(dt)
 + \wt Z\int_{T} h(t) \Phi_q\(\(\int_{T} h(v)\, m(dv)\)^{-1}
h(t)\) \, m(dt).\nn
\eea

It is easy to check thatÊ $x\,\Phi_q(x/\bb)$, or
equivalently, $x\,\Phi_q(x )$, is a convex
function for all $0\le q< \ff$. Consequently,
 it follows fromÊ Jensen's
inequalityÊ that
\bea
&&
\int_{T} h(t)\Phi_q\(\(\int_{T} h(v)\, m(dv)\)^{-1} h(t)\)
\, m(dt) \\
&&\qquad\quad \ge \int_{T} h(t)\, m(dt)\Phi_q(1)\nn.
\eea
Using this in (\ref{eee}) yields the inequality
\bas
&&\int_{T} |X(t)|h(t) \, m(dt) \\
&& \qquad \le D_q \wt Z
\int_{T} h(t) \Phi_q\(\(\int_{T} h(v)\, m(dv)\)^{-1} h(t)\)
\, m(dt).
\eas
where $D_q=1+(1/\Phi_q(1))$. Changing
$D_q \wt Z$ÊÊ to $Z$ gives \re{hX}.Ê \qed

\medskip\noindent
{\bf Proof of Theorem \ref{maj}} Using Lemma \ref{a2} it is easy toÊ complete
the proof of Theorem
\ref{maj} by following the proof ofÊÊ \cite[Theorem
5.2.6]{Fernique} orÊ Ê Ê \cite[TheoremÊ 6.3.3]{book}. We make some
commentsÊ regardingÊ Ê the proof inÊ \cite[TheoremÊ 6.3.3]{book}.
InÊÊ place of (6.73) we have  that for some $\al<\ff$
\bea
 E|X(t)-M_k(t)|&\le& \frac1{\mu_k(t)}\int_{B(t,D2^{-k})}
 E|X(t)-X(u)|\,\mu(du) \\
&\le& \frac1{\mu_k(t)}\int_{B(t,D2^{-k})}
 \al\|X(t)-X(u)\|_{\psi_q}\,\mu(du) \nn\\
&\le& \frac1{\mu_k(t)}\int_{B(t,D2^{-k})}
\al d(u,t)\,\mu(du)\le \al D2^{-k}\nn,
\eea
which is all we need to proceed with the proof. This follows because by
Jensen's Inequality, for any convex function $\Psi$,  
\be
E\Psi\({\al |X(t)-X(u)|\overÊ E|X(t)-X(u)|}\)\ge \Psi\(\al \).
\ee
LetÊ $\Psi(x)=\psi_q(x)$. Therefore, when $\psi_q(\al )\ge1$
 \be
 E|X(t)-X(u)|\le \al \|X(t)-X(u)\|_{\psi_q}.
\ee
It is easy to see that we can take $\al=1$ when $1\le q\le \ff$. When $q<1$
the reader can check that it suffices to take $\al=x_0$.

Ê When $1\le
q\le
\ff$Ê the rest of the adaptation of the proof of TheoremÊ 6.3.3 in \cite{book} is
completely apparent. When $0<q<1$ one gets as far as the expression on the
bottom of page 261 but with the measures multiplied by $G_q$, (and a
different constant following $Z$).Ê We need only be concerned if
$G_q<1$,Ê In this case we proceed as in \cite[(6.85)]{book} and note that
\be
\log\(1+\frac{x}{G_q}\)\le {\log\(1+(2/G_q)\)\over \log 2}\log x\qquad x\ge 2.
\ee
 Using this the proof can be completed.\qed

\bl \label{lem-5.4}For $0<q<1$, there exists a constantÊ $G_q>0$
such that
\be
\label{Yoq}
xy \le \psi_q(x) +Ê y\(2\log (1+y/G_q)\)^{1/q}ÊÊ \quad x, y \ge 0
\ee

\el

\ProofÊ It is easy to see that for all $p>0$ there exists a constantÊ $D_p>0$
for which
\be
\frac{e^s}{s^p }\ge D_p\(e^{s/2}-1\)\qquad \forall s\in R_+.\label{m1.4}
\ee
Taking $s=x^q$ this shows that there exists a
constantÊÊÊ $G_q>0$ such that
\be
{ \exp( x^q)\over x^{1-q}}\ge {G_q \over q} \( \exp (x^q/2) -1 \)\qquad
\forall x\in R_+.
\ee
By (\ref{psi-1})
\be
\psi'_q(x)=Ê {q\exp( x^q)\over x^{1-q}}\qquadÊ x> x_0.
\ee
Consequently
\be
\psi_q'(x)\geÊ G_qÊ \( \exp (x^q/2)-1 \)\qquad x>
x_0.\label{m1.2w}
\ee
Let $\La_{q}(y)$ be the right continuous inverse of $\psi'_q(x)$.
By (\ref{psi-1}) we have
$\La_{q}(y)=0$ for $y< K_q$ and $\La_{q}(y)=x_0$
for
$K_q\leÊ y\le D_{+}\psi_q(x_{0})
=Ê q\exp( x_{0}^q)/x_{0}^{1-q}$, the right hand derivative of
$\psi_q(x)$ at $x_{0}$.
In addition, by (\ref{m1.2w}) we see that
\be
\La_{q}(y)\leÊÊ \(2\log(1+y/G_q)\)^{1/q}\qquad y>{q\exp( x_{0}^q)
\over x_{0}^{1-q}}\label{m1.2wej}.
\ee
Therefore,ÊÊ decreasing $G_q$ if necessary, we have that
\be
\La_{q}(y)\leÊÊ \(2\log(1+y/G_q)\)^{1/q}\qquad \forall\, y\in
R_+,\label{m1.2we}
\ee
from which we get (\ref{Yoq}) by Young's Inequality and the obvious fact that
$\int_{0}^{y}\La_{q}(s)\,ds\leq y \La_{q}(y)$ since $\La_{q}(s)$ is
non-decreasing.
\qed

\section{Proof of Theorem \ref{theo-9.2}}\label{sec6} 

 Consider the Gaussian chaos ${\cal{X}}=\{:X_{h}^{k}:(g), h\in (0,1]\}$ defined in (\ref{chi}). It is clear that this process is continuous on $(0,1]$. Therefore, to show that it is continuous on $[0,1]$ it suffices to show that it is continuous on $[0,h_{0}]$ for some   $0<h_0<<1$. For
$h,h'\in[0,h_0]$ set
\begin{equation}
d(h,h'):=\|:X_{h}^{ k}:( g)-:X_{h'}^{ k}:( g)\|_{2}.\label{nor14.6m}
\end{equation}
  It follows from (\ref{nor4.6}) that
\be
\lim_{h,h'\to 0}d(h,h')=0.
\ee
Therefore, by \cite[Theorem 3.2.10]{Ginepena}
\be
\lim_{h,h'\to 0}\|:X_{h}^{ k}:( g)-:X_{h'}^{ k}:( g)\|_{\psi_{2/k}}=0.
\ee
Furthermore, the same theorem states that the $L^2$ and $L_{\psi_{2/k}}$
are equivalent.
Consequently
\begin{equation}
\lim_{h\rar 0}  :X_{h}^{ k}:( g)= :X_{0}^{ k}:( g)\qquad \mbox{in
$\, L_{\psi_{2/k}}$} \label{nor14.6ss}
\end{equation}
and
\be
\|:X_{h}^{ k}:( g)-:X_{h'}^{ k}:( g)\|_{\psi_{2/k}}\le C d(h,h')
\ee
for all h$,h'\in[0,1]$. 

We use Theorem \ref{maj} to   show that $\cal{X}$ is continuous on $(
[0,h_0],d)$. To do this we need estimates for $d$. We get one estimate from Lemma \ref{calc}. Th next lemma gives another estimate  for $d$.

\bl Under the hypotheses of Theorem \ref{theo-9.2},
 for any $h,h'>0$
\begin{equation}
d(h,h')=\|:X_{h}^{ k}:( g)-:X_{h'}^{ k}:( g)\|_{2}
\leq C\({|h-h'| \over hh'}\)^{1/2}.\label{nor14.7}
\end{equation}
\el

\Proof  Note that by (\ref{defin})
\begin{equation}
    X_{h}  ( g) \,\stackrel{def}{=}Ê :X_{h}^{ 1} ( g):\,\,  =\int  \(\frac{G(x+h)-G(x)}{h}
\) g(x)\,dx.
   \end{equation}
In addition it  is not hard to see that it follows from the definition of $\rho$ in (\ref{rho}),   that for $x'\leq x$, and $y'\leq y$
\begin{equation}
E\(G(x)-G(x')\)\(G(y)-G(y')\)=\int_{x'}^{x}\int_{y'}^{y}    \rho (t-s)   \,ds\,dt.\label{new.2star}
\end{equation}
(Details are  given  in \cite[Lemma 2.2]{ae}.) Therefore
\begin{eqnarray}
\lefteqn{ E\( X_{ h} ( g)X_{ h'}(\wt g)\)
\label{new.4}}\\
&&= {1 \over h}{1 \over h'}\int\int E\(\(G(x+h)-G(x)\)\(G(y+h')-G(y)\)\)\,g( x)\,dx\,\wt g( y)\,dy
\nonumber\\
&&  = \int\int\lc {1 \over h}\int_{x}^{x+h}{1 \over h'}\int_{y}^{y+h'}    \rho (t-s)   \,ds\,dt\rc\,g( x)
\,\wt g( y)\,dx\,dy.\nonumber 
\end{eqnarray}
 Let   $(X,Y)$ be a
two dimensional Gaussian random variable.     
By   \cite[Theorem 3.9]{Janson}  
\begin{equation} E(:X^{ k}::Y^{  j}:)=k!(E( XY))^{ k}\de_{
k,j}.\label{18.3a}
\end{equation} 
Using this and (\ref{new.4}) we see that   
\begin{eqnarray}\lefteqn{
E (:X_{ h}^{k}:( g):X^{k}_{ h'}:( \wt g))
\label{nor4.3caa}}\\
&& = k!\int \int \({1 \over h}\int_{ x}^{ x+h}
{1 \over h'}\int_{ y}^{ y+h'}\rho( s-t)\,dt \,ds\)^{k}\,g( x)\wt g( y)\,dx\,dy.
\nonumber
\end{eqnarray}

Set
\bea 
B_{z}(h,h')&=&{1 \over h}\int_{0}^{ h}{1 \over h'}\int_{ 0}^{ h'}\rho(z+ s-t)\,dt
\,ds\label{bzdef}\\
&=&{ \si^2(z+h)+ \si^2(z-h')- \si^2(z+h-h')-\si^2(z) 
\over 2hh'}.  \nn
\eea
By (\ref{nor4.3caa}) and a change of variables   we
have
 \begin{eqnarray}
\lefteqn{\|:X_{h}^{ k}:( g)-:X_{h'}^{ k}:( g)\|_{2}^{2}
 \label{nor14.7f}}\\
 &&=\int\int \lc \(B_{z}(h,h)\)^{k}
-\(B_{z}(h,h')\)^{k}-\(B_{z}(h',h)\)^{k}+\(B_{z}(h',h')\)^{k}\rc \nn\\
&&\hspace{3.7in}g(x) g(
y)\,dx\,dy  .
\nonumber
 \end{eqnarray}

In addition
\begin{eqnarray}
\lefteqn{B_{z}(h,h)-B_{z}(h,h')
\label{nor14.21r}}\\
&&=\(\frac{1}{h^2}-\frac{1}{hh'}\) h^2B_{z}(h,h)+\frac{1}{hh'}\(
h^2B_{z}(h,h)-hh'B_{z}(h,h')\) \nn\\ 
&&  = {1 \over 2hh'}\(\si^2(z+h-h')+\si^2(z-h)-\si^2(z-h')-\si^2(z)\)
 \nonumber\\
&& \hspace{.5 in}+\({1 \over 2h^2}-{1 \over 2hh'}\)\( \si^2(z+h)+\si^2(z-h) -2
\si^2(z)\). \nonumber
\end{eqnarray}
We write this as \begin{eqnarray}
\lefteqn{B_{z}(h,h)-B_{z}(h,h')
\label{nor14.21r8}}\\
&&  = {1 \over2 hh'}\(\si^2(z+h-h')-\si^2(z)+\si^2(z-h)-\si^2(z-h')\)
 \nonumber\\
&& \hspace{.5 in}  +{(h-h') \over 2h'h^{2}}\(2\si^2(z)-\si^2(z+h)-\si^2(z-h)\).
\nonumber
\end{eqnarray}
Since $\si^2$ and $(\si^2)'$ are   bounded  we need only use the mean
value theorem,   on four differences,  to see that   for $h,h'>0$
\begin{equation}
|B_{z}(h,h)-B_{z}(h,h')|\leq C{|h-h'| \over h'h}.\label{nor14.21r9}
\end{equation}

Note that
\begin{equation}
B^{k}_{z}(h,h)-B^{k}_{z}(h,h')=\sum_{j=0}^{k-1}B^{j}_{z}(h,h)
(B_{z}(h,h)-B_{z}(h,h'))B^{
k-j-1}_{z}(h,h')
\label{nor14.7ew}.
\end{equation}
Therefore
\begin{equation}
|B^{k}_{z}(h,h)-B^{k}_{z}(h,h')|\le  C{|h-h'| \over h'h}
\sum_{j=0}^{k-1}B^{j}_{z}(h,h)B^{k-j-1}_{z}(h,h')
\label{nor14.7ej}.
\end{equation}

  			Let $f_{h}(x)={1 \over h}1_{[0,h]}(x)$ so that
 the first line in the definition (\ref{bzdef}) of $B_{z}(h,h')$ can be written as
\begin{equation}
B_{z}(h,h')=\int \int\rho(z+ s-t)f_{h}(s)f_{h'}(t)\,ds
\,dt.\label{14.2.1h}
\end{equation}
Using Fubini's Theorem we see that
\begin{eqnarray}
&&  \int\int B^{j}_{z}(h,h)B^{k-j-1}_{z}(h,h')\, g(x) g(
y)\,dx\,dy
\label{3.7ah}\\
&& \qquad   =\int\!\!\int\(\int\dots\int\) \prod_{i=1}^{k-1} \rho (x+v_i-y-w_i)
\prod_{i=1}^{j} f_{h}(v_i) f_{h}(w_i)\,dv_i\,dw_i \nn\\ 
 &&\hspace{1.5in}\prod_{i=j+1}^{k-1}
f_{h}(v_i) f_{h'}(w_i)\,dv_i\,dw_i \,\, g(x)g(y)\,dx\,dy\nn\\&&
\qquad =\int\dots\int\(\int\!\!\int \prod_{i=1}^{k-1}
\rho(x-y+v_i-w_i) g(x)g(y)\,dx\,dy\)\nn\\
&&\hspace{1in} \prod_{i=1}^{j}
f_{h}(v_i) f_{h}(w_i)\,dv_i\,dw_i\prod_{i=j+1}^{k-1}
f_{h}(v_i) f_{h'}(w_i)\,dv_i\,dw_i\nn\\
&&\qquad \leq C\nn
\end{eqnarray}
  where  $C $ is a finite constant that is independent of $h$ and $h'$.
In the last step we use  the generalized Holder's inequality and  the fact  that
$\rho\in L^{k}_{loc}$ and $g\in\BB_{0}(R_{+)}$, to get
\begin{equation}
\int\!\!\int \prod_{i=1}^{k-1}
\rho(x-y+v_i-w_i) g(x)g(y)\,dx\,dy\leq C.\label{3.7ai}
\end{equation}
Using (\ref{3.7ah}) together with (\ref{nor14.7ej}) we obtain
\be
 \int\int \bigg| \(B_{z}(h,h)\)^{k}
-\(B_{z}(h,h')\)^{k}  \bigg|\, g(x) g(
y)\,dx\,dy \le C' {|h-h'| \over h'h}.\label{14.2.1}
\ee
Clearly the integral of the other two terms in (\ref{nor14.7f}) has the same
bound. Thus we get  (\ref{nor14.7}). \qed

\medskip\noindent{\bf Proof of Theorem \ref{theo-9.2} }
It follows from (\ref{.03})   that for any $h>0$
\begin{equation}
d(h,0)=\|:X_{h}^{ k}:( g)-:X_{0}^{ k}:( g)\|_{2}\leq Ch^{\de/2}\label{nor14.7a}
\end{equation} (The constant $C$ actually depends on $k$, but we take  $k$ fixed.)
We use this bound as well as the one in (\ref{nor14.7}).

By Theorem \ref{maj}  to prove that
${\cal X}$ is continuous it suffices to show that
\be
\sup_{h\in[0,h_0]}\int_0^{Kh_0^{\de/2}}\(\log{1\over
\la(B_d(h,u))}\)^{k/2}\,du<\ff,\label{mm.1f}
\ee
 and
\be
\lim_{\ep\to 0}\sup_{h\in[0,h_0]}\int_0^\ep\(\log{1\over
\la(B_d(h,u))}\)^{k/2}\,du=0,\label{mm.1}
\ee
where $\la$ is Lebesgue measure. (Theorem \ref{maj} requires a probability
measure. Rather than bothering to renormalize we need only observe that its
conclusions also hold for positive measures with mass less than one.)  

We  now verify (\ref{mm.1f}) and (\ref{mm.1}).Ê We pick $h_0$ so that
$Kh_0^{\de/2}$ is very small.Ê LetÊ $h\in(0,h_0]$. Note that by (\ref{nor14.7})
we have $h'\in B_d(h,u)$ when
$|(h-h')/(hh')|^{1/2}\le u/C$, or, equivalently, whenÊ $| h-h' |\le
hh'(u/C)^{2}$.Ê (We take $C\ge 1$.) Since $h'\leq 1$ and $u\le Kh_0^{\de/2} $,
we see thatÊ ÊÊ on
$\{h':| h-h' |\le hh'u^{2}\}$ weÊ have
$h'>h/2$. Therefore
\begin{equation}
B_d(h,u)\supseteq \{h':| h-h' |\le h^{2}u^{2}/(2C^2)\},\label{mm.2a}
\end{equation}
so that theÊ Lebesgue measure of
$B_d(h,u)$ is at least $h^{2}u^{2}/(2C^2)$. Consequently
  for any $h\in (0,h_0]$ and $u\le Kh_0^{\de/2} $
\be
\log{1\over
\la(B_d(h,u))}\le2\( \log{1\over
h}+\log{1\over
u}+\log C\).\label{92.7}
\ee
Therefore for any $h\in (0,h_0]$ and  $w\le Kh_0^{\de/2} $
\bea
&&\int_{0}^{w}\(\log{1\over
\la(B_d(h,u))}\) ^{k/2}\,du\label{mm.2wj}\\
&&\qquad \leÊ ÊÊ C'Ê w\( \log{1\over
h}+\log{1\over w}+\log C\)^{k/2}\leq
C''Ê w\( \log{1\over
h}+\log{1\over w}\)^{k/2}\nn.
\eea

 LetÊ $h\in (0,h_0]$ and $v\le Kh_0^{\de/2} $
and suppose thatÊ
$h^{\de/4}\ge v$. Then by (\ref{mm.2wj}) and   the monotonicity of $
\log 1/h$
\bea
&&\int_{0}^{v}\(\log{1\over
\la(B_d(h,u))}\) ^{k/2}\,du \leÊÊ Cv\( \log{1\over
v}\)^{k/2}.\label{mm.2w}
\eea
(The constants are not necessarily the same at each stage.)
Now suppose that
$h^{\de/4}< v$. In this case using (\ref{mm.2wj})Ê Ê with $w=h^{\de/4}$ we
have
\bea
&&\int^{h^{\de/4}}_0\(\log{1\over
\la(B_d(h,u))}\) ^{k/2}\,du \label{mm.5}\\
&&\qquad \le C_\de h^{\de/4}\(\log 1/h\)^{k/2} \nn\\
&&\qquad \le C'_\de h^{\de/4}\(\log 1/h^{\de/4}\)^{k/2} \leÊÊ K_\de v\(
\log{1\over v}\)^{k/2}\nn.
\eea
 (Here we use the monotonicity of   $x\(\log 1/x\)^{k/2}$.)

Now consider
\be
\int_{ h^{\de/4}}^{v}\(\log{1\over
\la(B_d(h,u))}\) ^{k/2}\,du.
\ee
Since $d(h,x)<d(h,0)+d(x,0)$, we see by  
(\ref{nor14.7a}) that
$\{x\in B_d(h,u) \}$ when
$Ch^{\de/2}+Cx^{\de/2}\le u$, or, equivalently, when $x\le
C'(u-Ch^{\de/2})^{2/\de}$. Since $u\ge h^{\de/4}$, we see that for small $h$,
(which we can always achieve by taking $Kh_0^{\de/2}$ sufficiently small) we
have
$x\le
C'(u-Ch^{\de/2})^{2/\de}$ whenever $x\le C''u ^{2/\de}$, for some $C''>0$.
Consequently $\la(B_d(h,u))\ge Ku ^{2/\de}$ and
\be
\int_{ h^{\de/4}}^{v}\(\log{1\over
\la(B_d(h,u))}\) ^{k/2}\,du\leÊ C v\(\log 1/v\)^{k/2} \label{mm.3s}.
\ee
Using (\ref{nor14.7a}) it is elementary toÊ see that
\be
\int_{0} ^{v}\(\log{1\over
\la(B_d(0,u))}\) ^{k/2}\,du \le
 C v\(\log 1/v\)^{k/2}.\label{mm.4}
\ee
CombiningÊ (\ref{mm.2w}), (\ref{mm.5}), (\ref{mm.3s})Ê and (\ref{mm.4}) we get that for any
$v\le Kh_0^{\de/2} $
\be
\sup_{h\in[0,h_0]}\int_0^v\(\log{1\over
\la(B_d(h,u))}\)^{k/2}\,du\leqÊ C v\(\log 1/v\)^{k/2}.\label{mm.1jx}
\ee
 The statements in (\ref{mm.1f}) and
 (\ref{mm.1}) follow immediately.\qedÊ

\def\noopsort#1{} \def\printfirst#1#2{#1}
\def\singleletter#1{#1}
             \def\switchargs#1#2{#2#1}
\def\bibsameauth{\leavevmode\vrule height .1ex
             depth 0pt width 2.3em\relax\,}
\makeatletter
\renewcommand{\@biblabel}[1]{\hfill#1.}\makeatother
\newcommand{\bysame}{\leavevmode\hbox to3em{\hrulefill}\,}

\end{document}